\documentclass [11pt]{article}
\begin{document}
\thispagestyle{empty}
\null\vspace{-1cm}
\medskip
\vspace{1.75cm}
\begin{center}

           {\bf Score sets in oriented bipartite graphs}

\vspace {1in}

S. Pirzada$^{1}$, T. A. Naikoo$^{2}$ and T. A. Chishti$^{3}$\\
\bigskip

 $^{1,2}$Department of Mathematics, University of Kashmir, Srinagar-190006, India.\\

$^{1}$Email :     sdpirzada@yahoo.co.in\\

$^{2}$Email :      tariqnaikoo@rediffmail.com\\

$^{3}$Mathematics Section, Centre of Distance Education, University of Kashmir, Srinagar-190006, India.\\

$^{3}$Email :      chishtita@yahoo.co.in\\ 
\end{center}

\vspace {1in}

{\bf Abstract}. The set A of distinct scores of the vertices of an oriented bipartite graph D(U, V) is called its score set. We consider the following question: given a finite, nonempty set A of positive integers, is there an oriented bipartite graph D(U, V) such that score set of D(U, V) is A?  We conjecture that there is an affirmative answer, and verify this conjecture when $\mid A \mid$  = 1, 2, 3, or when A is a geometric or arithmetic progression.\\ 
 \bigskip
\vfill

\newpage

\vspace{0.15in}

{\bf 1. Introduction}. An oriented graph is a digraph with no symmetric pairs of directed arcs and without loops. Let D be an oriented graph with vertex set V = $\{v_{1}, v_{2},~.~.~.~, v_{n}\}$, and let $d_{v}^+$  and $d_{v}^-$  denote the outdegree and indegree respectively of a vertex v. Avery [1] defined $a_{v} = n - 1 + d_{v}^+ - d_{v}^-$, the score of v, so that $0 \leq a_{v} \leq 2n-2$. Then, the sequence  $[a_{1}, a_{2},~.~.~.~,  a_{n}]$ in non-decreasing order is called the score sequence of D.\\

\indent Avery [1] obtained the following criterion for score sequences in oriented graphs .\\

{\bf Theorem 1.1.} A non-decreasing sequence of non-negative integers $[a_{1}, a_{2}, ~.~.~.~ , a_{n}]$ is the score sequence of an oriented graph if and only if\\
\begin{center}
$\sum_{i=1}^{k}a_{i}$ $\geq ~ k(k-1)$, ~	for $1 \leq k\leq n$,
\end{center}
with equality when k = n.\\

\indent Pirzada and Naikoo [7] obtained the following results for score sets in oriented graphs.\\

{\bf Theorem 1.2.}  Let A = $\{a, ad, ad^{2}, ~.~.~.~ , ad^{n}\}$, where a and d are positive integers with a $>$ 0 and d $>$ 1. Then, there exists an oriented graph D with score set A, except for a = 1, d = 2,  n $>$ 0 and for a = 1, d = 3,  n $>$ 0.\\

{\bf Theorem 1.3.}  If $a_{1}, a_{2}, ~.~.~.~ , a_{n}$ are n non-negative integers with $a_{1} < a_{2} <~.~.~.~< a_{n}$, then there exists an oriented graph D with score set A = $\{a^{'}_{1}, a{'}_{2}, ~.~.~.~ , a{'}_{n}\}$, where\\

\begin{center} $a^{'}_{i}= \bigg\{^{a_{i-1}+a_{i}+1,~~~~for~i~>~1,}_{a_{i},~~~~for~i~=~1.}$\\ 
\end{center}

\indent The study of score sets in tournaments (complete oriented graphs) can be found in [2, 5, 8, 10, 11].\\
\indent  An oriented bipartite graph is the result of assigning a direction to each edge of a simple bipartite graph. Let U = $\{u_{1}, u_{2} , ~.~.~.~ , u_{m}\}$ and V = $\{v_{1}, v_{2} , ~.~.~.~ , v_{n}\}$ be the parts of an oriented bipartite graph D(U, V). For any vertex x in   D(U, V) , let $d_{x}^{+}$ and $d_{x}^{-}$ respectively be the outdegree and indegree of x. Define $a_{u}=n+d_{u}^{+}-d_{u}^{-}$ and $b_{v}=m+d_{v}^{+}-d_{v}^{-}$  respectively as the scores of u in U and v in V. Clearly , $0 \leq a_{u} \leq 2n$ and $0 \leq b_{v} \leq 2m$. The sequences $[a_{1}, a_{2} , ~.~.~.~ , a_{m}]~~ and~~ [b_{1}, b_{2} , ~.~.~.~ , b_{n}]$ in non-decreasing order are called the score sequences of  D(U, V).\\

 The following result due to Pirzada, Merajuddin and Yin [4] is the bipartite version of Theorem 1.1.\\

{\bf Theorem 1.4.} Two non-decreasing sequences $[a_{1} , a_{2} ,~.~.~.~, a_{m}]~~ and~~ [b_{1}, b_{2},~.~.~.~, b_{n}]$ of non-negative integers are the score sequences of some oriented bipartite graph if and only if\\
\begin{center}
$\sum_{i=1}^{p} a_{i}$ + $\sum_{j=1}^{q} b_{j}\geq 2pq$,~ for $1 \leq p \leq m$ and $1 \leq q \leq n$,
\end{center} 
with equality when p = m and q = n.\\

          The study of score sets for bipartite tournaments (complete oriented bipartite graphs) can be found in [3, 9, 12] and for k-partite tournaments (complete oriented k-partite graphs) in [6].\\

{\bf 2. Score sets in oriented bipartite graphs}\\

{\bf Definition.} The set A of distinct scores of the vertices in an oriented bipartite graph D(U, V) is called its score set. 
            If there is an arc from a vertex u to a vertex v, then we say that vertex u dominates vertex v. \\
 
               We have the following results.\\

{\bf Theorem 2.1.}  Every singleton or doubleton set of positive integers is a score set of some oriented bipartite graph.\\
{\bf Proof. Case I.}  Let  A = \{a\}, where a is a positive integer. When a is even, construct an oriented bipartite graph D(U, V) as follows.\\
\indent Let 
\begin{center}
U = $X_{1} \cup  X_{2}$,\\

V = $Y_{1} \cup Y_{2}$
\end{center}    
with $X_{1}  \cap X_{2}  = \phi, Y_{1}  \cap Y_{2}  = \phi,  | X_{1} | = | X_{2} | = | Y_{1} | = | Y_{2} | =  \frac {a}{2}$. Let every vertex of $X_{i}$ dominates each vertex of $Y_{i}$, and every vertex of $Y_{i}$ dominates each vertex of $X_{j}$ whenever $i \neq j$ so that we get the oriented bipartite graph D(U, V) with \\
\begin{center}
$| U | = | X_{1} | + | X_{2} | = | Y_{1} | + | Y_{2} | = | V | =\frac{a}{2}+\frac{a}{2}$ = a,\\
\end{center}
and the scores of vertices \\
 \indent~~ $a_{x_{1}}= | V | + | Y_{1} | - | Y_{2} | = | U | + | X_{1} | - | X_{2} | = a_{y_{2}} = a + \frac{a}{2} -\frac{a}{2}  = a$,  for all $x_{1} \in X_{1}, y_{2} \in Y_{2}$\\
and ~~$a_{x_{2}}  = | V | + | Y_{2} | - | Y_{1} | = | U | + | X_{2} | - | X_{1} | = a_{y_{1}} = a +  \frac{a}{2} -\frac{a}{2}  = a$,  for all $x_{2} \in X_{2}, y_{1} \in Y_{1}$.\\
Therefore, score set of  D(U, V) is A = \{a\}.\\

               Now, when a is odd, construct an oriented bipartite graph D(U, V) as follows.\\ \indent  Let
\begin{center}
                      $U = X_{1} \cup X_{2} \cup \{x\}$,\\  
                          $V = Y_{1} \cup Y_{2} \cup \{y\}$\\
\end{center}    
with $X_{1}  \cap X_{2}  = \phi, X_{i}  \cap \{x\} = \phi, Y_{1}  \cap Y_{2}  = \phi, Y_{i}  \cap \{y\}  = \phi,  | X_{1} | = | X_{2} | =   | Y_{1} | = | Y_{2} | =\frac{a-1}{2}$. Let every vertex of $X_{i}$ dominates each vertex of $Y_{i}$, and every vertex of $Y_{i}$ dominates each vertex of $X_{j}$ whenever i $\neq$ j so that we get the oriented bipartite graph D(U, V) with\\
\begin{center} 
$| U | = | X_{1} | + | X_{2} | + | \{x\} |  = | Y_{1} | + | Y_{2} | + | \{y\} | = | V | = \frac{a-1}{2} +\frac{a-1}{2}  + 1 = a$,\\
\end{center}
and the scores of vertices \\
 \indent ~~ $ a_{x_{1}}    = | V | + | Y_{1} | - | Y_{2} | = | U | + | X_{1} | - | X_{2} | =a_{y_{2}}$ = $a +\frac{a-1}{2}  -\frac{a-1}{2}  = a$,  for all $x_{1} \in X_{1}, y_{2} \in Y_{2},$\\
 \indent ~~ $ a_{x_{2}}   = | V | + | Y_{2} | - | Y_{1} | = | U | + | X_{2} | - | X_{1} | = a_{y_{1}} = a + \frac{a-1}{2}  -\frac{a-1}{2} = a$,  for all $x_{2} \in X_{2}, y_{1} \in Y_{1}$\\
and     $a_{x} = | V | + 0 - 0 = | U | + 0 - 0 = a_{y} = a$, for the vertices x and y.\\
Thus, score set of  D(U, V) is A = \{a\}.\\

            ~~~Note that an empty oriented bipartite graph D(U, V) with $| U | = | V | = a$ has also score set A = \{a\}.\\
  
{\bf Case II.}  Let  $A = \{a_{1}, a_{2}\}$, where  $a_{1}$ and $a_{2}$ are positive integers with $a_{1} <  a_{2}$. As in case I, there exists an oriented bipartite graph D(U, V) with $| U | = | V | = a_{1}$, and the scores of vertices
              $a_{u} = a_{v} = a_{1}$,  for all $u \in U, v \in V$. \\
\indent Since $a_{2} > a_{1}$ or $a_{2} - a_{1} > 0$, construct oriented bipartite graph $D(U_{1}, V_{1})$ as follows.\\
\indent  Let  $U_{1} = U \cup X,  V_{1} = V,  U \cap X = \phi,  | X | = a_{2} - a_{1}$. Let there be no arc between the vertices of V and X, so that we get  the oriented bipartite graph D$(U_{1}, V_{1})$ with\\
\begin{center}
                                $| U_{1}| = | U | + | X | = a_{1} + a_{2} - a_{1} = a_{2},  | V_{1}| = a_{1}$,\\
\end{center}
and the scores of vertices\\
\indent ~       $a_{u} = a_{1},$  for all $u \in U$,\\
\indent ~      $a_{x} = | V_{1}| + 0 - 0 = a_{1}$,  for all $x \in X$,\\
 and          $a_{v} = a_{1} + | X | = a_{1} + a_{2} - a_{1} = a_{2},$   for all $v \in V$.\\
Hence, score set of  $D(U_{1}, V_{1} )$ is  A =$\{a_{1}, a_{2}\}$.\\  
 \indent  Again, note that an empty oriented bipartite graph D(U, V) with $| U | = a_{1}, | V | = a_{2}$ has also score set $A = \{a_{1} , a_{2}\}$.\\
 
\indent {\bf Theorem 2.2.}   Every  set  of  three  positive integers  is  a  score  set  of  some oriented bipartite graph.\\
{\bf Proof.}  Let $A = \{a_{1}, a_{2}, a_{3}\}$, where $a_{1}, a_{2}, a_{3}$ are positive integers with $a_{1} < a_{2} < a_{3}$.\\
\indent  First assume $a_{3} > 2a_{2}$  so that  $a_{3} - 2a_{2} > 0$, and since $a_{2} > a_{1}$, therefore $a_{3} - 2a_{1} > 0$. Now, construct an oriented bipartite graph D(U, V) as follows.\\
                Let   $U =  X_{1} \cup X_{2},  V = Y_{1} \cup Y_{2}$  with $X_{1} \cap X_{2} = \phi ,  Y_{1} \cap Y_{2} = \phi$,  $| X_{1} | = a_{2} , | X_{2}  | = a_{3} - 2a_{2}, | Y_{1} | = a_{1}, | Y_{2} | = a_{3} - 2a_{1}$. Let every vertex of $X_{2}$ dominates each vertex of $Y_{1}$, and every vertex of $Y_{2}$ dominates each vertex of $X_{1}$, so that we get the oriented bipartite graph D(U, V) with\\
\begin{center}
                                 $| U | = | X_{1} | + | X_{2} | = a_{2} + a_{3} - 2a_{2} = a_{3} - a_{2}$,\\ 
                                 $| V | = | Y_{1} | + | Y_{2} | = a_{1} + a_{3} - 2a_{1} = a_{3} - a_{1}$,\\ 
\end{center}
and the scores of vertices\\	
\indent ~ $a_{x_{1}}=\mid V \mid +0-(a_{3}-2a_{1})=a_{3}-a_{1}-a_{3}+2a_{1}=a_{1},$~~~for~all~ $x_{1} \in X_{1}$,\\ \indent ~ $a_{x_{2}}=\mid V \mid +a_{1}-0=a_{3}-a_{1}+a_{1}=a_{3},$~~~for all~ $x_{2} \in X_{2},$\\ \indent ~ $a_{y_{1}}=\mid U \mid +0-(a_{3}-2a_{2})=a_{3}-a_{2}-a_{3}+2a_{2}=a_{2},$~~~for all~ $y_{1} \in Y_{1}$,\\ and $a_{y_{2}}=\mid U \mid +a_{2}-0=a_{3}-a_{2}+a_{2}=a_{3},$~~~for all~ $y_{2} \in Y_{2}.$\\
Therefore, score set of D(U, V) is $A = \{a_{1} , a_{2} , a_{3}\}$.\\
\indent  Now, assume $a_{3} \leq 2a_{2}$ so that $2a_{2} - a_{3} \geq 0$. Construct an oriented bipartite graph D(U, V) as follows.\\
\indent  Let  $U = X_{1}, V = Y_{1} \cup Y_{2}$  with $Y_{1} \cap Y_{2} = \phi, | X_{1} | = a_{2}, | Y_{1} | = a_{1},     | Y_{2} | = a_{2} - a_{1}$. Let every vertex of $Y_{2}$ dominates $a_{3} - a_{2}$ vertices of $X_{1}$ (out of $a_{2}$), so that we get the oriented bipartite graph D(U, V) with\\
\begin{center}
        $| U | = | X_{1} | = a_{2},  | V | = | Y_{1} | + | Y_{2} | = a_{1} + a_{2} -a_{1} = a_{2}$,
\end{center}
and the scores of vertices\\ \indent ~ $a_{x_{1}}=\mid V \mid +0-(a_{2}-a_{1})=a_{2}-a_{2}+a_{1}= a_{1}$,~~~for the $a_{3}-a_{2}$ vertices of $X_{1}$,\\ \indent ~ $a_{x^{'}_{1}}=\mid V \mid +0-0=a_{2},$ for the remaining $a_{2}-(a_{3}-a_{2})=2a_{2}-a_{3}$  vertices of $X_{1}$,\\ \indent ~ $a_{y_{1}}=\mid U \mid +0-0=a_{2},$~~~for~all~$y_{1} \in Y_{1},$\\ and~$a_{y_{2}}=\mid U \mid +a_{3}-a_{2}-0=a_{2}+a_{3}-a_{2}=a_{3},$~~~for~all~$y_{2} \in Y_{2}.$\\              
Thus, score set of D(U, V) is $A = \{a_{1}, a_{2}, a_{3}\}$.\\
        
\indent The next result shows that every set of positive integers in geometric progression is a score set of some oriented bipartite graph.\\

{\bf Theorem 2.3.}  Let $A = \{a, ad, ad^{2}, ~.~.~.~ , ad^{n}\}$, where a and d are positive integers with $a > 0$ and $d > 1$. Then, there exists an oriented bipartite graph  with score set A.\\ 
{\bf Proof.}  First assume $d > 2$. Induct on n. If n = 0, then by Theorem 2.1, there exists an oriented bipartite graph D(U, V) with score set A = \{a\}.\\
\indent For n = 1, construct an oriented bipartite graph D(U, V) as follows.\\ 
\indent Let  U = $X_{1} \cup X_{2}$,  V = $Y_{1} \cup Y_{2}$ with $X_{1} \cap X_{2} = \phi, Y_{1} \cap Y_{2} = \phi, | X_{1} | = | Y_{1} | = a, | X_{2} | = | Y_{2} | = ad - 2a > 0$ as  $a > 0,~ d > 2$. Let every vertex of $X_{2}$ dominates each vertex of $Y_{1}$, and every vertex of $Y_{2}$ dominates each vertex of  $X_{1}$, so that we get the oriented bipartite graph D(U, V) with\\
\begin{center}
                                $| U | = | X_{1} | + | X_{2} | = a + ad - 2a = ad - a$,\\ 
                                $| V | = | Y_{1} | + | Y_{2} | = a + ad - 2a = ad - a$,\\
\end{center}
and the scores of vertices \\ \indent ~ $ a_{x_{1}} = | V | + 0 - (ad-2a) = ad-a-ad+2a = a$,  for all $x_{1} \in X_{1},$\\ \indent ~ $ a_{x_{2}}   =  | V | + a-0 = ad-a+a = ad,$ for all $x_{2} \in X_{2}$,\\ \indent ~ $a_{y_{1}}  = | U | + 0 - (ad-2a) = ad-a-ad+2a = a$,  for all $y_{1} \in Y_{1},$\\ and $ a_{y_{2}}   =  | U | + a-0 = ad-a+a = ad,$ for all $y_{2} \in Y_{2}$.\\  Thus, score set of D(U, V) is $A = \{a, ad\}$.\\
\indent Assume the result to be true for all $p \geq 1$. We show that the result is true for p + 1.\\ \indent Let a and d be positive integers with $a > 0$ and $d > 2$. Therefore, by induction hypothesis, there exists an oriented bipartite graph D(U, V) with \\
\begin {center}
$\mid U \mid = \mid V \mid = ad^{p}-(ad^{p-1}-ad^{p-2}+ ~.~.~.~ (-1)^{p+1}a)$,
\end {center}                             
and  $a , ad , ad^{2} , ~.~.~.~ , ad^{p}$ as the scores of the vertices of D(U, V). As $a > 0 , ~   d > 2$ , therefore $ad^{p+1}-2(ad^{p}-(ad^{p-1}-ad^{p-2}+ ~.~.~.~ (-1)^{p+1}a))>0$.   Now, construct an oriented bipartite graph $D(U_{1} ,~ V_{1} )$ as follows.\\
\indent  Let  $U_{1} = U \cup X ,  V_{1} = V \cup Y$ with $U \cap X = \phi , V \cap Y = \phi$,\\
\begin {center}
$\mid X \mid = \mid Y \mid = ad^{p+1}-2(ad^{p}-(ad^{p-1}-ad^{p-2}+ ~.~.~.~ (-1)^{p+1}a))$.\\
\end {center}                   
Let every vertex of X dominates each vertex of V, and every vertex of Y dominates each vertex of  U, so that we get the oriented bipartite graph $D(U_{1} , V_{1})$ with\\$| U_{1} | = | U | + | X |=| V | + | Y | = | V_{1} |$\\ = $ad^{p}-(ad^{p-1}-ad^{p-2}+ ~.~.~.~ (-1)^{p+1}a) + ad^{p+1}-2(ad^{p}-(ad^{p-1}-ad^{p-2}+ ~.~.~.~ (-1)^{p+1}a))$\\ = $ad^{p+1}-(ad^{p}-(ad^{p-1}-ad^{p-2}+ ~.~.~.~ (-1)^{p+1}a))$,\\ 
and since $| X | = | Y |$, therefore  $a + | X | - | X | = a,~ ad + | X | - | X | = ad ,~ ad^{2} +        | X | - | X | = ad^{2}, ~.~.~.~ , ad^{p} + | X | - | X | = ad^{p}$ are the scores of the vertices of U and V, and \\ \indent ~ $ a_{x}   =  | V_{1} | + | V | - 0 = | U_{1} | + | U | - 0 = a_{y} = ad^{p+1}-(ad^{p}-(ad^{p-1}-ad^{p-2}+ ~.~.~.~ (-1)^{p+1}a))+ ad^{p}-(ad^{p-1}-ad^{p-2}+ ~.~.~.~ (-1)^{p-1}a)= ad^{p+1}$, for all $x \in X, y \in Y.$\\ 
Therefore, score set of $D(U_{1},~V_{1} )$ is $A = \{a , ad , ad^{2} , ~.~.~.~ , ad^{p} , ad^{p+1}\}$.\\
\indent Now, assume d = 2. Then the set A becomes  $A = \{a, 2a, 2^{2}a, ~.~.~.~ , 2^{n}a \}$. Construct an oriented bipartite graph D(U, V) as follows. \\
\indent Let 
\begin {center}
   $U = X_{0} \cup X_{1} \cup X_{3} \cup X_{4} \cup ~.~.~.~ \cup X_{n}$,\\
               $V = Y_{0} \cup Y_{2} \cup Y_{3} \cup Y_{4} \cup ~.~.~.~ \cup Y_{n}$\\
\end {center}	
with $X_{i}  \cap X_{j}  = \phi$, $Y_{i}  \cap Y_{j}  = \phi$ $(i \neq j)$. Let $| X_{0} | = | X_{1} | = | Y_{0} | = | Y_{2} | = a$, \\and for $3 \leq i \leq n $\\
\begin {center}
$| X_{i} | = | Y_{i} | = 2^{i}a-2(\sum_{j=0,j \neq 2}^{i-1}| X_{j} |),~~~~~~~~~~~~~~~~~~~~~~~~~~~~~~~(2.3.1)$\\
\end {center}                                                     
which is clearly greater than zero. Let every vertex of $X_{i}$ dominates each vertex of $Y_{j}$ whenever $i > j$, and every vertex of $Y_{i}$ dominates each vertex of $X_{j}$ whenever $i > j$, so that we get the oriented bipartite graph D(U, V) with the scores of vertices\\ \indent ~$a_{x_{0}} = | V |+0- \sum_{j=2}^{n}| Y_{j} | = \sum_{j=0,j \neq 1}^{n}| Y_{j} |- \sum_{j=2}^{n}| Y_{j} | =| Y_{0} |=a,$  for all $ x_{0} \in X_{0}$,\\ \indent ~ $a_{x_{1}} = | V |+| Y_{0} |- \sum_{j=2}^{n}| Y_{j} | = \sum_{j=0,j \neq 1}^{n}| Y_{j} |+a- \sum_{j=2}^{n}| Y_{j} | =| Y_{0} |+a=2a,$  for all $ x_{1} \in X_{1}$,\\ \indent ~ $a_{y_{0}} = | U |+0- \sum_{j=1,j \neq 2}^{n}| X_{j} |=\sum_{j=0,j \neq 2}^{n}| X_{j} |-\sum_{j=1,j \neq 2}^{n}| X_{j} |=| X_{0} | =a,$  for all $ y_{0} \in Y_{0}$,\\ \indent ~ $a_{y_{2}} = | U |+| X_{0} |+| X_{1} |- \sum_{j=3}^{n}| X_{j} | = \sum_{j=0,j \neq 2}^{n}| X_{j} |+a+a- \sum_{j=3}^{n}| X_{j} | =| X_{0} |+| X_{1} |+2a=a+a+2a=4a,$  for all $ y_{2} \in Y_{2}$,\\              
and for $3 \leq i \leq n$\\ \indent ~ $a_{x_{i}} = | V |+\sum_{j=0,j \neq 1}^{i-1}| Y_{j} | - \sum_{j=i+1}^{n}| Y_{j} |=| U |+ \sum_{j=0,j \neq 2}^{i-1}| X_{j} |- \sum_{j=i+1}^{n}| X_{j} | =a_{y_{i}}= \sum_{j=0,j \neq 2}^{n} \mid X_{j} \mid +   \sum_{j=0,j \neq 2}^{i-1} \mid X_{j} \mid -\sum _{j=i+1}^{n} \mid X_{j} \mid = \sum_{j=0,j \neq 2}^{i} \mid X_{j} \mid +   \sum_{j=0,j \neq 2}^{i-1} \mid X_{j} \mid = 2\sum_{j=0,j \neq 2}^{i-1} \mid X_{j} \mid + \mid X_{i} \mid = 2\sum_{j=0,j \neq 2}^{i-1} \mid X_{j} \mid + 2^{i}a - 2(\sum_{j=0,j \neq 2}^{i-1} \mid X_{j} \mid)$\\ \indent ~~~~~~~~~~~~~~~~~~~~~~~~~~~~~~~~~~~~~~~~~~~~~~~~~~~~~ (By equation (2.3.1))\\ \indent ~ ~~~~ $= 2^{i}a$, for all $x_{i} \in X_{i}, y_{i} \in Y_{i}$.\\    
Therefore , score set of D(U, V) is  $A = \{a , 2a , 2^{2}a , ~.~.~.~ , 2^{n}a\}$.\\
  
\indent The next result shows that every set of positive integers in arithmetic progression is a score set for some oriented bipartite graph.\\

\indent {\bf Theorem 2.4}  Let $A = \{ a , a + d , a + 2d , ~.~.~.~ , a + nd \}$, where a and d are positive integers. Then, there exists an oriented bipartite graph with score  set A.\\
{\bf Proof.(a).}  Let d $>$ a  so that  $d - a > 0$. Construct an oriented bipartite graph D(U, V) as follows. \\
\indent   Let  
\begin {center} 
   $U = X_{0} \cup X_{1} \cup  ~.~.~.~ \cup X_{n}$,\\
              $V = Y_{0} \cup Y_{1} \cup ~.~.~.~ \cup Y_{n}$\\ 
\end {center}
with $X_{i}  \cap X_{j}  = \phi , Y_{i}  \cap Y_{j}  = \phi ( i \neq j )$, and for $0 \leq i \leq n$\\
\begin{center}
 $ \mid X_{i} \mid =  \mid Y_{i} \mid = \bigg\{^{a,~~~~~~~~if ~i ~is ~even,}_{d-a,~~~~~~if ~i ~is~ odd.}~~~~~~~~~~~~~~~~~~~~~~~~~~~~~~~~~(2.4.1)$ 
\end{center}                                                 
Let every vertex of $X_{i}$ dominates each vertex of $Y_{j}$ whenever $i > j$, and every vertex of $Y_{i}$ dominates each vertex of $X_{j}$ whenever $i > j$, so that we get the oriented bipartite graph D(U, V) with\\
\indent $ \mid U \mid = \sum _{i=0}^{n} \mid X_{i} \mid = \sum _{i=0}^{n} \mid Y_{i} \mid = \mid V \mid$\\ \indent ~~~~ = $\bigg\{^{a+d-a+a+d-a+~.~.~.~+d-a+a,~~~~~~~~if~n~is~even,}_{a+d-a+a+d-a+~.~.~.~+a+d-a,~~~~~~~~if~n~is~odd,}$\\ \indent ~~~~ = $\bigg\{^{(\frac {n}{2}+1)a+ \frac {n}{2}(d-a),~~~~~~~~if~n~is~even,}_{(\frac {n+1}{2})a+ (\frac {n+1}{2})(d-a),~~~~~if~n~is~odd,}$\\ \indent ~~~~ = $ \bigg\{^ {\frac {nd}{2}+a,~~~~~~~~if~n~is~even,}_{( \frac {n+1}{2})d,~~~~~~~~if~n~is~odd,}~~~~~~~~~~~~~~~~~~~~~~~~~~~~~~~~~~~~~~~~~~~~~~~~~~~~~(2.4.2)$\\         
and the scores of vertices\\ \indent $a_{x_{0}} = | V |+0- \sum_{j=1}^{n}| Y_{j} | = | U |+0-\sum_{j=1}^{n}| X_{j} |= a_{y_{0}} = \sum_{j=0}^{n}| Y_{j} | -\sum_{j=1}^{n}=| Y_{j} |= \mid Y_{0} \mid =a,$  for all $ x_{0} \in X_{0},~y_{0} \in Y_{0}$,\\          
and for $1 \leq i \leq n$\\ \indent $a_{x_{i}} = | V |+\sum_{j=0}^{i-1}| Y_{j} | - \sum_{j=i+1}^{n}| Y_{j} |=| U |+ \sum_{j=0}^{i-1}| X_{j} |- \sum_{j=i+1}^{n}| X_{j} |= a_{y_{i}} = \sum_{j=0}^{n}| X_{j} | + \sum_{j=0}^{i-1}| X_{j} |-\sum_{j=i+1}^{n}| X_{j} |= \sum_{j=0}^{i}| X_{j} |+ \sum_{j=0}^{i-1}| X_{j} |=2\sum_{j=0}^{i-1}| X_{j} |+| X_{i} |$\\ \indent ~~~ = $\bigg\{^{2\sum_{j=0}^{i-1}| X_{j} |+a,~~~~~~~~~~~~if~i~is~even,}_{2\sum_{j=0}^{i-1}| X_{j} |+d-a,~~~~~~~~~if~i~is~odd,}$~~~~~~~~~~(By equation (2.4.1))\\ \indent ~~~ =$ \bigg\{^{2(\frac {i-1+1}{2})d+a,~~~~~~~~~~~~if~i~is~even,}_{2((\frac {i-1}{2})d+a)+d-a,~~~~~~~if~i~is~odd,}$~~~~~~~~~~~~~(By equation (2.4.2))\\ \indent ~~~ = $\bigg\{^{a+id,~~~~~~~~if~i~is~even,}_{a+id,~~~~~~~~if~i~is~odd.}$\\ 
That is, $a_{x_{i}}=a_{y_{i}}=a+id$, for all $x_{i} \in X_{i}$,  $y_{i} \in Y_{i}$ where $1 \leq i \leq n$. Therefore, score set of D(U, V) is  A = \{a, a + d, a + 2d ,~.~.~.~  , a + nd\}.\\
{\bf (b).} Let d = a. Then the set A becomes A = \{a, 2a, 3a,~.~.~.~ , (n + 1)a\}. For n = 0, the result follows from Theorem 2.1 . Now, assume $n \geq 1$.\\
\indent If n is odd, say $n = 2k - 1$ where $k \geq 1$, then construct an oriented bipartite graph D(U, V) as follows.\\
 \indent Let   
\begin {center}
$U = X_{0} \cup X_{1} \cup X_{3} \cup ~.~.~.~ \cup X_{2k-3} \cup X_{2k-1}$,\\ $V = Y_{0} \cup Y_{2} \cup Y_{4} \cup ~.~.~.~ \cup Y_{2k-2}$\\
\end {center}
with $X_{i}  \cap X_{j}  = \phi , Y_{i}  \cap Y_{j}  = \phi ( i \neq j ) , ~and~ | X_{i} | = | Y_{j} | = a , ~for~ all                   ~i \in \{0 , 1 , 3 , ~.~.~.~ , 2k-1\} ,  j \in \{ 0 , 2 , 4 , ~.~.~.~ , 2k-2 \}$ . Let every vertex of $X_{i}$ dominates each vertex of $Y_{j}$ whenever $i > j$, and every vertex of $Y_{i}$ dominates each vertex of $X_{j}$ whenever $i > j > 0$, so that we get the oriented bipartite graph D(U, V) with\\
\begin {center}
 $\mid U \mid = \sum _{j \in \{0 , 1 , 3 , ~.~.~.~ , 2k-1\}} \mid X_{j} \mid = a+(\frac {2k-1+1}{2})a=a+ka,$\\$ \mid V \mid = \sum _{j \in \{0 , 2 , 4 , ~.~.~.~ , 2k-2\}} \mid Y_{j} \mid = a+(\frac {2k-2}{2})a=ka,$\\
\end {center}                             
and the scores of vertices\\ \indent ~~ $a_{x_{0}}=\mid V \mid +0-0=ka,$  for all $x_{0} \in X_{0},$\\for  $i \in \{ 1 , 3 , ~.~.~.~ , 2k -1 \}$\\ \indent ~~ $a_{x_{i}} = | V |+| Y_{0} |+\sum_{j \in \{2 , 4 , ~.~.~.~ , i-1\}}| Y_{j} | - \sum_{j \in \{i+1,i+3 , ~.~.~.~ , 2k-2\}}| Y_{j} |=ka+a+(\frac {i-1}{2})a-(\frac {2k-2-(i-1)}{2})a = ka+a+ia-a-ka+a=(i+1)a,$  for all $x_{i}\in X_{i}$,\\ \indent ~~ $a_{y_{0}} = | U |+0-\sum_{j \in \{1, 3, ~.~.~.~ , 2k-1\}}| X_{j} | = a+ka-(\frac {2k-1+1}{2})a=a,$  for all $y_{0} \in Y_{0},$\\            
and for  $i \in \{ 2 , 4 , ~.~.~.~ , 2k - 2 \}$\\ \indent ~~ $a_{y_{i}} = | U |+\sum_{j \in \{1 , 3 , ~.~.~.~ , i-1\}}| X_{j} | - \sum_{j \in \{i+1,i+3 , ~.~.~.~ , 2k-1\}}| X_{j} |=a+ka+(\frac {i-1+1}{2})a-(\frac {2k-1+1-(i-1+1)}{2})a =a+ka+\frac {ia}{2}-ka+ \frac {ia}{2}=(i+1)a,~~for~all~y_{i} \in Y_{i}.$\\         
Thus, score set of D(U, V ) is  $A = \{ a , 2a , 3a , ~.~.~.~ , (2k - 1) a , 2ka \}$.\\
\indent Now, if n is even, say n = 2k where $k \geq 1$, then construct an oriented bipartite graph D(U, V) as follows.\\ 
\indent Let   
\begin {center}
  $U = X_{0} \cup X_{1} \cup X_{3} \cup ~.~.~.~ \cup X_{2k-1}$,\\
             $V = Y_{0} \cup Y_{2} \cup Y_{4} \cup ~.~.~.~ \cup Y_{2k}$\\
\end {center} 
with $X_{i}  \cap X_{j}  = \phi , Y_{i}  \cap Y_{j}  = \phi ~~( i \neq j )$ , and $| X_{i} | = | Y_{j} | = a$ , for all  $i \in \{0 , 1 , 3 , ~.~.~.~ , 2k-1\} ,~  j \in \{ 0 , 2 , 4 , ~.~.~.~ , 2k \}$. Let every vertex of $X_{i}$ dominates each vertex of $Y_{j}$ whenever $i > j$, and every vertex of $Y_{i}$ dominates each vertex of $X_{j }$ whenever $i > j > 0$, so that we get the oriented bipartite graph D(U, V) with ( as in above )  $| U | = a + ka ,~  | V | = ka + a = a + ka$, and the scores of vertices\\
\indent ~~  $a_{x_{0}} = ka + | Y_{2k} | = ka + a = ( k +1 ) a,$  for all $x_{0} \in X_{0}$,\\
for $i \in \{ 1 , 3 , ~.~.~.~ , 2k - 1 \}$\\ 
\indent ~~  $a_{x_{i}} = ( i + 1 ) a,$  for all $x_{i} \in X_{i}$,\\
\indent ~~  $a_{y_{0}}= a,$  for all $y_{0} \in Y_{0}$,\\
for $i \in \{ 2 , 4 , ~.~.~.~ , 2k - 2 \}$\\
\indent ~~  $a_{y_{i}} = ( i + 1 ) a ,  ~~for~ all~ y_{i} \in Y_{i}$,\\
and $a_{y_{2k}} = | U |+\sum_{j \in \{1, 3, ~.~.~.~ , 2k-1\}}| X_{j} |-0 = a+ka+(\frac {2k-1+1}{2})a=(2k+1)a,$  for all $y_{2k} \in Y_{2k}.$\\               
Thus, score set of D(U, V ) is  $A = \{ a , 2a , 3a , ~.~.~.~ , 2ka , ( 2k +1) a \}$.\\
{\bf (c).}  Let $d < a$ so that $a - d > 0$ . For n = 0 or 1, the result follows from Theorem 2.1 . Now, assume that $n \geq 2$.\\
\indent If n is even, say n = 2k where $k \geq 1$, then construct an oriented bipartite graph D(U, V) as follows.\\ 
\indent  Let     
\begin {center}
$U = X_{0} \cup X_{1} \cup X_{3} \cup ~.~.~.~ \cup X_{2k-1}$,\\
$V = Y_{0} \cup Y_{2} \cup Y_{4} \cup ~.~.~.~ \cup Y_{2k}$\\
\end {center}
with $X_{i} \cap X_{j}  = \phi , ~ Y_{i}  \cap Y_{j}  = \phi ~~( i \neq j ),| X_{0} | = | Y_{0} | = a,$ and $| X_{i} | = | Y_{j} | = d,$ for all $i \in \{1 , 3 , ~.~.~.~ , 2k-1\}, ~ j \in \{ 2 , 4 , ~.~.~.~ , 2k \}$. Let every vertex of $X_{i}$ dominates each vertex of $Y_{j}$ whenever $i > j > 1$, every vertex of $X_{i}$ dominates d vertices of $Y_{0}$ ( out of a ) whenever $i > 2$, and every vertex of $Y_{i}$ dominates each vertex of $X_{j}$ whenever $i > j > 0$, so that we get the oriented bipartite graph D(U, V) with\\
\begin {center}
$ \mid U \mid = \sum _{j \in \{0 , 1 , 3 , ~.~.~.~ , 2k-1\}} \mid X_{j} \mid = a+(\frac {2k-1+1}{2})d=a+kd,$\\$ \mid V \mid = \sum _{j \in \{0 , 2 , 4 , ~.~.~.~ , 2k\}} \mid Y_{j} \mid = a+(\frac {2k}{2})d=a+kd,$\\
\end {center}            
and the scores of vertices\\ \indent ~~ $a_{x_{0}}=\mid V \mid +0-0=a+kd,$  for all $x_{0} \in X_{0},$\\ \indent ~~  $a_{x_{1}}=\mid V \mid +0- \sum _{j \in \{2 , 4 , ~.~.~.~ , 2k\}} \mid Y_{j} \mid = a+kd-(\frac {2k}{2})d=a,$  for all $x_{1} \in X_{1},$\\for  $i \in \{3 , 5, ~.~.~.~ , 2k -1 \}$\\ \indent ~~  $a_{x_{i}} = | V |+d+\sum_{j \in \{2 , 4 , ~.~.~.~ , i-1\}}| Y_{j} | - \sum_{j \in \{i+1,i+3 , ~.~.~.~ , 2k\}}| Y_{j} |=a+kd+d+(\frac {i-1}{2})d-(\frac {2k-(i-1)}{2})d = a+kd+d+(i-1)d-kd=a+id,~~for~all~x_{i}\in X_{i}$,\\ \indent ~~  $a_{y_{0}} = | U |+0-0=a+kd,$  for the $a-d$ vertices of $Y_{0}$,\\ \indent ~~  $a_{y^{'}_{0}}= \mid U \mid + 0 -\sum_{j \in \{3 , 5 , ~.~.~.~ , 2k-1\}}| X_{j} | = a+kd-(\frac {2k-1+1-(1+1)}{2})d=a+kd-kd+d=a+d,$  for the remaining d vertices of $Y_{0},$\\ and for  $i \in \{2 , 4, ~.~.~.~ , 2k \}$\\ \indent ~~  $a_{y_{i}} = | U |+\sum_{j \in \{1 , 3 , ~.~.~.~ , i-1\}}| X_{j} | - \sum_{j \in \{i+1,i+3 , ~.~.~.~ , 2k-1\}}| X_{j} |=a+kd+(\frac {i-1+1}{2})d-(\frac {2k-1+1-(i-1+1)}{2})d = a+kd+\frac {id}{2}-kd+\frac{id}{2}=a+id,$  for all $y_{i} \in Y_{i}.$\\Therefore, score set of  D(U, V ) is  $A = \{a , a + d , a + 2d , ~.~.~.~ , a + ( 2k -1) d , a + 2kd \}$.\\ \indent Now, if n is odd, say n = 2k + 1 where $k \geq 1$, then construct an oriented bipartite graph D(U, V ) as follows.\\ 
\indent   Let 
\begin {center}
$U = X_{0} \cup X_{1} \cup X_{3} \cup ~.~.~.~ \cup X_{2k-1} \cup X_{2k+1}$,\\
               $V = Y_{0} \cup Y_{2} \cup Y_{4} \cup ~.~.~.~ \cup Y_{2k}$\\
\end {center}
with $X_{i} \cap X_{j}  = \phi , Y_{i}  \cap Y_{j}  = \phi ~( i \neq j ) , | X_{0} | = | Y_{0} | = a$, and $| X_{i} | = | Y_{j} | = d$, for all $i \in \{ 1 , 3 , ~.~.~.~ , 2k +1\} ,  j \in \{ 2 , 4 , ~.~.~. , 2k \}$. Let every vertex of $X_{i}$ dominates each vertex of $Y_{j}$ whenever $i > j > 1$, every vertex of $X_{i}$ dominates d vertices of $Y_{0}$ ( out of a ) whenever $i > 2$, and every vertex of $Y_{i}$ dominates each vertex of $X_{j}$ whenever $i > j > 0$, so that we get the oriented bipartite graph D(U, V) with (as in above)  $| U | = a + kd +d = a + ( k + 1 ) d , | V | = a + kd$, and the scores of vertices\\ \indent ~~ $a_{x_{0}}=a+kd,$  for all $x_{0} \in X_{0},$\\ \indent ~~ $a_{x_{1}}=a,$  for all $x_{1} \in X_{1},$\\for  $i \in \{3 , 5, ~.~.~.~ , 2k-1 \}$\\ \indent ~~ $a_{x_{i}} = a+id,$  for all $x_{i} \in X_{i},$\\ \indent ~~ $a_{x_{2k+1}}=| V |+d+\sum_{j \in \{2 , 4 , ~.~.~.~ , 2k\}}| Y_{j} | -0=a+kd+d+(\frac {2k}{2})d=a+(2k+1)d,$  for all $x_{2k+1} \in X_{2k+1},$\\ \indent ~~ $a_{y_{0}} = a+kd+| X_{2k+1} |=a+kd+d=a+(k+1)d,$  for the $a-d$ vertices of $Y_{0}$,\\$a_{y^{'}_{0}}= a+d,$  for the remaining d vertices of $Y_{0}$,\\and for $i \in \{2, 4,~.~.~.~,2k\}$\\ \indent ~~ $a_{y_{i}}=a+id,$  for all $y_{i} \in Y_{i}.$\\    
Hence, score set of  D(U, V ) is  $A = \{a , a + d , a + 2d , ~.~.~.~ , a + 2kd , a + (2k + 1) d \}$, and the proof is complete.\\

                         ¦
{\bf Remark.} We note that Theorems 2.1, 2.2, and 2.4 cannot be extended to state that any set of nonnegative integers  A  is a score set of some oriented bipartite graph when  $| A |$ = 1, 2, 3, or when A is an arithmetic progression, for instance, there is no oriented bipartite graph with score set  \{0\}, \{0, 1\},  or \{0, 1, 2\}.\\

            We conclude with the following conjecture.\\

\indent {\bf Conjecture.} Every finite set of positive integers is a score set for some oriented bipartite graph.\\

\bigskip
{\bf References}\\

\bigskip
\indent [1]  Avery, P., Score sequences of oriented graphs, J. Graph Theory, Vol. 15, No. 3 (1991)251- \indent ~~~~257.\\ 
\indent [2]  Hager, M., On score sets for tournaments, J. Discrete Mathematics 58 (1986) 25-34.\\
\indent [3]  Petrovic, V., On bipartite score sets, Univ. u Novom Sadu Zb. Rad. Prirod. Mat. Fak. \indent ~~~~ Ser. Mat. 13 (1983) 297-303.\\ 
\indent [4]  Pirzada, S., Merajuddin, Yin,  J., On the scores of oriented bipartite graphs, J.\\  
\indent ~~~~Mathematical Study, Vol. 33, No. 4 (2000) 354 - 359.\\ 
\indent [5]  Pirzada, S., Naikoo, T. A., On score sets in tournaments, Vietnam J. of\\ \indent ~~~~Mathematics, Vol. 34 (2006) To appear.\\ 
\indent [6]  Pirzada, S., Naikoo, T. A., Score sets in k-partite tournaments, J. of Applied\\ \indent ~~~~  Mathematics and Computing (2006), To appear.\\ 
\indent [7]  Pirzada, S., Naikoo, T. A., Score sets in oriented graphs, To appear.\\ 
\indent [8]  Reid, K. B., Score sets for tournaments, Congressus Numerantium XXI, Proceedings of \\ \indent ~~~~the Ninth Southeastern Conference on Combinatorics, Graph Theory, and Computing \\ \indent ~~~~(1978) 607-618.\\   
\indent [9]  Wayland, K., Bipartite score sets, Canadian Mathematical Bulletin, Vol. 26, No. 3 \\ \indent ~~~~(1983)~273-279.\\ 
\indent [10] Yao, T. X., Reid's conjecture on score sets in tournaments (in Chinese), Kexue \\ \indent ~~~~~Tongbao 33 (1988) 481-484.\\ 
\indent [11] Yao, T. X., On Reid's conjecture of score sets for tournaments, Chinese Sci. Bull. 34 \\ \indent ~~~~~(1989) 804-808.\\ 
\indent [12] Yao, T. X., Score sets of bipartite tournaments, Nanjing Daxue Xucbao Ziran Kexue  \\ \indent ~~~~~Ban 26  (1990) 19-23.\\

\end{document}